\documentclass[a4paper,oneside,11pt]{article}
\usepackage{verbatim,color}
\usepackage{textcomp}
\usepackage[english]{babel}
\usepackage{amssymb}
\usepackage{amsmath}

\definecolor{Red}{cmyk}{0,1,1,0}
\definecolor{Blue}{cmyk}{1,1,0,0}

\renewcommand{\leq}{\leqslant}
\renewcommand{\geq}{\geqslant}

\usepackage{amsmath,amsfonts,amssymb,amsthm,mathrsfs,stackrel}

\def\build#1_#2^#3{\mathrel{
\mathop{\kern 0pt#1}\limits_{#2}^{#3}}}

\newcommand{\E}{\mathbb{E}}
\renewcommand{\a}{\alpha}
\renewcommand{\t}{\tau}
\renewcommand{\l}{\lambda}

\newcommand{\ZZ}{\mathbb{Z}}
\theoremstyle{plain}
\newtheorem{theorem}{Theorem}

\newtheorem{lemma}{Lemma}

\theoremstyle{definition}

\theoremstyle{remark}

\begin{document}
\title{A remark on monotonicity in Bernoulli bond Percolation}

\author{Bernardo N.B. de Lima\footnote{Departamento de Matem{\'a}tica, Universidade Federal de Minas Gerais, Av. Ant\^onio
Carlos 6627 C.P. 702 CEP 30123-970 Belo Horizonte-MG, Brazil}, Aldo Procacci$^{*}$ and
R\'emy Sanchis$^{*}$}

\maketitle
\begin{abstract}
Consider an anisotropic
independent bond percolation model on the $d$-dimensional hypercubic lattice, $d\geq 2$, with parameter $p$. We show that the two point connectivity function $P_{p}(\{(0,\dots,0)\leftrightarrow (n,0,\dots,0)\})$ is a monotone function in $n$ when the parameter $p$ is close enough to 0. Analogously, we show that truncated connectivity function  $P_{p}(\{(0,\dots,0)\leftrightarrow (n,0,\dots,0), (0,\dots,0)\nleftrightarrow\infty\})$ is also a monotone function in $n$ when $p$ is close to 1.
\end{abstract}

{\footnotesize


}

\section{Introduction and main result}

Consider an ordinary Bernoulli bond percolation model, with
parameter $p\in[0,1]$ on the graph $\mathbb{L}^d=(\ZZ^d,\E)$ for $d\ge 2$, with
$\E=\{e=\langle x, y\rangle: x,y\in \ZZ^d,~|x-y|_1=1\}$. That is
each bond is {\em open} independently with probability $p$,
otherwise it is {\em closed} with probability $1-p$. Thus, this
model is described by the probability space $(\Omega, {\cal F},
P_p)$ where $\Omega=\{0,1\}^{\E}$, ${\cal F}$ is the
$\sigma$-algebra generated by the cylinder sets in $\Omega$ and
$P_{p}=\prod_{e\in \E} \mu(e) $ is the product of Bernoulli measures
with parameter $p$.

Given two vertices $x,y\in\ZZ^d$, we use the standard notation $(x\leftrightarrow y)$ to denote the set of configurations $\omega\in\Omega$
such that $x$ is connected to $y$ by a path of open bonds. Given the parameter $p\in[0,1]$ and $n\in\mathbb{N}$, we define
the two-point connectivity function $\tau_p(n)=:P_p((0,\dots,0)\leftrightarrow (n,0,\dots,0))$. It is still
an open question to prove if $\tau_p(n)$ is monotone in $n\in\mathbb{N}$ for all values of $p$;
this problem was told to one of us (B.N.B.L.) by J. van den Berg \cite{vdB1}.
We also define the truncated two-point function $\tau^f_p(n)=:P_p((0,\dots,0)\leftrightarrow (n,0,\dots,0),(0,\dots,0)\nleftrightarrow\infty)$
as the probability of the set of configurations where the origin is connected by open paths to the vertex $(n,0,\dots,0)$
but the origin is to connected by open paths to at most finitely many vertices.

From now on, as usual $p_c=p_c(d)$ is the percolation threshold for ordinary Bernoulli percolation on $\mathbb{L}^d$. We remind that $\tau^f_p(n)$ is the interesting quantity in the supercritical phase since $\tau_p(n)$ does not decay at all when  $p>p_c$.

The main result of this note is the following theorem:
\begin{theorem} \label{prin} Consider the Bernoulli bond percolation on $ \mathbb{L}^d$ with $d\ge 2$, then:\\
i) There exists $p^\prime>0$ (depending upon the dimension $d$), such that $$\tau_p(n)>\tau_p(n+1),\ \forall n\in\mathbb{N},\ \forall p<p^\prime.$$
ii) There exists $p^{\prime\prime}<1$ (depending upon the dimension $d$), such that $$\tau^f_p(n)>\tau^f_p(n+1),\ \forall n\in\mathbb{N},\ \forall p>p^{\prime\prime}.$$
\end{theorem}

An analogous question was conjectured by Hammersley and Welsh \cite{HW} for the monotonicity of expected passage times in the context of first-passage percolation, and this question is also still open. There are some partial results like \cite{Go}, \cite{H}, \cite{AW} and \cite{A}. One very interesting negative result is due van den Berg \cite{vdB2}, he consider first-passage percolation on the graph $\ZZ_+\times\ZZ$ and proves that the expect passage time from the origin to $(2,0)$ is less than the expect passage time from the origin to $(1,0)$.

Another similar inequality in the context  of oriented percolation was obtained by E. Andjel and M. Sued in \cite{AS}.

\section{Proof of Theorem \ref{prin}}

Essentially the proof combines well established  estimates of the Ornstein-Zernike behavior for the correlation functions for large values of $n$
and  classical estimates via polymer expansion for small values of $n$. Here, as usual $\xi(p)$  and $\xi^f(p)$ are the correlation length and the truncated correlation length defined as (see equations (6.54) and (8.56) in \cite{Gr}):
$$\xi(p)=\left[\lim_{n\rightarrow\infty}-\frac{\log \tau_p(n)}{n}\right]^{-1}$$
and
$$\xi^f(p)=\left[\lim_{n\rightarrow\infty}-\frac{\log \tau^f_p(n)}{n}\right]^{-1}.$$

\subsection{Proof of i)}
The proof is based on two results. The first one concerns the Ornstein-Zernike decay of the
two-point connectivity function $\t_p(n)$ in the whole subcritical phase and was originally proved by and Campanino, Chayes and Chayes  \cite{CCC}.
The second result concerning upper and lower bounds for the  two-point connectivity function $\t_p(n)$ in the highly subcritical phase,
was obtained in \cite{PSBS} via polymer expansion. Hereafter  all constants depend on $d$.

\begin{lemma}\label{ccc}[Theorem 6.2 (II) of \cite{CCC}] Consider independent bond percolation on $\mathbb{L}^d,\ d\geq 2$. For all
$p<p_c$ there exists $\a(p)>0$, $K_2(p)>0$ such that

$$ \t_p(n)=\frac{K_2(p)}{\left(\a(p)\pi
n\right)^{\frac{d-1}{2}}}\exp\left\{-\frac{n}{\xi
(p)}\right\}\left[1+O(n^{-1})\right],\ \forall n\in\mathbb{N} .$$
\end{lemma}

\begin{lemma}\label{psbs}[Eq. (6.3) of \cite{PSBS}] Consider independent bond percolation on
$\ZZ^d$, $d\geq 2$. There exist $p_0>0$  close enough to zero and constants $C_1,C_2>0$, such that $\forall p<p_0$, it holds:

$$p^n(1-p)^{C_1n}\leq \t_p(n) \leq p^n (1+C_2p)^{n/2},\ \forall n\in\mathbb{N}.$$
\end{lemma}
\noindent
By Lemma \ref{ccc}, it holds that
$$
\frac{\t_p(n)}{\t_p(n+1)}=(1+1/n)^{\frac{d-1}{2}}e^{\frac{1}{\xi(p)}}\left(\frac{1+O(n^{-1})}{1+O((n+1)^{-1})}\right),\ \forall n\in\mathbb{N}.
$$
\noindent
Observe that $-\frac{C(p)}{n}\leq [O(n^{-1})]\leq \frac{C(p)}{n}$,
where $C(p)$ is a bounded function of $p$, at least in some
interval $[0,p_1]$, for $p_1$ small. Therefore,
$$\frac{1+O(n^{-1})}{1+O((n+1)^{-1})}\geq \frac{1-C/n}{1+C/(n+1)}\geq
1-2C/n,$$ where $C=\sup\{C(p); 0\leq p\leq p_1\}$.

\noindent
Hence, there exists $n_0=n_0(p_1)$ such that, for
$n\geq n_0,$ we have
$$(1+1/n)^\frac{d-1}{2}(1-2C/n)>e^{-\frac{1}{\xi(p)}},\ \forall p\in [0,p_1].$$

\noindent
Therefore the monotonicity of $\t_p(n)$ is proved for all $n\geq
n_0$ and all $p\in[0,p_1]$.
Now we prove the monotonicity of $\t_p(n)$ for all $n\leq n_0$. By Lemma \ref{psbs}, for all $n\leq n_0$ and $p\in[0,p_0]$, it holds that

$$\frac{\t_p(n)}{\t_p(n+1)}\geq \frac{1}{p} \frac{(1-p)^{C_1 n}}{(1+C_2
p)^{(n+1)/2}}\geq \frac{1}{p} \left(\frac{(1-p)^{C_1}}{(1+C_2
p)^{1/2}}\right)^{n+1}\geq \frac{1}{p}
\left(\frac{(1-p)^{C_1}}{(1+C_2 p)^{1/2}}\right)^{n_0+1},$$
where $n_0$ is the positive integer defined above.

\noindent
Thus, there exists $p_2\in[0,p_0]$, small enough, such that for all $p<p_2$, we have
$$\frac{1}{p} \left(\frac{(1-p)^{C_1}}{(1+C_2
p)^{1/2}}\right)^{n_0+1}>1.$$

\noindent
We have thus proved the monotonicity of $\t_p(n)$ for all $n\leq
n_0$ in the interval $[0,p_2]$. Then, part {\em i)} of Theorem \ref{prin} is proved taking $p^\prime=\min\{p_1,p_2\}$.

\subsection{ Proof of ii)}

The proof for $d\geq 3$ is based on analogous results for the highly supercritical phase. The Ornstein-Zernike decay
of the truncated two-point function for the  Bernoulli bond percolation in the highly supercritical phase  (Lemma \ref{bps} below) was originally  obtained in \cite{BPS}, while the upper and lower bounds for the  truncated two-point function $\t_p(n)$ in the highly supercritical phase (Lemma \ref{psbs2} below) was once again given in \cite{PSBS}.

The proof of ii) $d\geq 3$ then follows  from the two lemmas below performing  the same steps of part {\em i)} with minor modifications. From now on, let
$\l=\l(p)=\frac{1-p}{p}$.

\begin{lemma}\label{bps}[Theorem 1.1 of \cite{BPS}] For $d\geq 3$,
there exists $p_3<1$, close enough to 1, such that for all $p\in[p_3,1)$, it holds that

$$\t_p^f(n)=\frac{(1-p)^{4d-d^2-1}}{(2\pi)\frac{d-1}{2}} (1+g(p)) \frac{e^{-\frac{n}{\xi^f(p)}}}{n^{\frac{d-1}{2}}}(1+O(n^{-1}))$$
where $g(p)$ is an analitic function for $p\in[p_3,1]$.
\end{lemma}

\begin{lemma}\label{psbs2}[Theorem 5.1 of \cite{PSBS}]
For $d\geq 3$, there exists a constant $C>0$ and $p_4<1$, close enough to 1, such that,
for all $p\in[p_4,1)$, it holds that

$$ \left(\frac{\l}{(1+\l)^2}\right)^{2(d-1)(n+1)+2}\leq t_p^f(n)\leq 2 (\l\sqrt{1+C\l})^{2(d-1)(n+1)+2}.$$
\end{lemma}

Finally, for the case $d=2$, we use two analogous results recently obtained. The first one, about Ornstein-Zernike behaviour, was given
in \cite{CIL} and it is stated below as Lemma \ref{cil}. The second one, on bounds upper and lower bonds for truncated two-point function
was obtained in \cite{CLS}, is stated below as Lemma \ref{cls}.

Using these  two lemmas here below, the proof of ii) for $d=2$ follows the same lines as in the subcritical case.

\begin{lemma}\label{cil}[Theorem 1.1 of \cite{CIL}] For $d=2$, there
are $\psi(p)>0$ and $p_5<1$, close to 1, such that
$$\t_p^f(n)=\psi(p)\frac{e^{-\frac{n}{\xi^f(p)}}}{n^2}(1+f(p,n))$$

where $-c(p)f(n)\leq f(p,n)\leq c(p)f(n)$, with $f(n)\rightarrow 0$
as $n\rightarrow\infty$ and $c(p)$ is bounded in $p\in [p_5,1]$.
\end{lemma}

\begin{lemma}\label{cls}[Proposition 2 of \cite{CLS}] For $\ZZ^2$, there exists $p_6<1$, close enough to 1.
such that, for all $p\in[p_6,1)$, it holds that
$$\l^{2n+2}p^{2n}\leq \t_p^f(n)\leq \l^{2n+2}\left[\frac{(4^3\l)^{n/2+1}}{1-4^3\l}+(1+12\l)^n\right].$$
\end{lemma}

\section*{Acknowledgments}
B.N.B.L. and A.P.  have been partially supported by the Brazilian  agencies
Conselho Nacional de Desenvolvimento Cient\'\i fico e Tecnol\'ogico
(CNPq) and  Funda\c c\~ao de Amparo \`a  Pesquisa do Estado de Minas Gerais (FAPEMIG - Programa de Pesquisador Mineiro).
R.S. has been partially supported by Conselho Nacional de Desenvolvimento Cient\'\i fico e Tecnol\'ogico
(CNPq).

\end{document}